\documentclass[11pt,a4paper,reqno]{amsart}
\usepackage[marginratio=1:1,totalwidth=15.75cm,totalheight=22.275cm]{geometry}
\usepackage[T1]{fontenc}
\usepackage{lmodern}
\usepackage[utf8]{inputenc}
\usepackage{amsmath,amssymb,amsthm,enumerate,amscd}
\usepackage[shortlabels]{enumitem}
\usepackage{graphicx}
\usepackage{csquotes}

\usepackage[dvipsnames]{xcolor}
\usepackage[colorlinks=true,linkcolor=blue!60!black,citecolor=teal!80!black,urlcolor=RoyalBlue]{hyperref}
\usepackage{tikz-cd}
\usepackage{amscd}
\usepackage{mathrsfs}

\numberwithin{equation}{section}

\usepackage{soul}
\usepackage{mathscinet}

\newcommand{\N}{\mathbb{N}}

\newcommand{\C}{\mathbb{C}}

\newcommand{\D}{\mathbb{D}}

\def\Chat{\widehat{\mathbb{C}}}

\def\C{{\mathbb{C}}}
\def\D{{\mathbb{D}}}
\def\N{{\mathbb{N}}}

\theoremstyle{plain}
\newtheorem{thm}{Theorem}[section]
\newtheorem{lemma}[thm]{Lemma}

\theoremstyle{definition}

\theoremstyle{remark}

\title[A Retrospective on Eremenko and Lyubich's entire functions]{From pathological to paradigmatic: A Retrospective on Eremenko and Lyubich's entire functions}

\author[N. Fagella]{Núria Fagella}
\author[L. Pardo-Simón]{Leticia Pardo-Simón}

\address{\noindent Dept. de Matemàtiques i Informàtica\\ Universitat de Barcelona\\ Catalonia\\ Spain\\
	\newline  Centre de Recerca Matemàtica\\ Bellaterra\\ Catalonia\\ Spain.
	\textsc{\newline \indent 
		\href{https://orcid.org/0000-0002-5466-0579%
		}{\includegraphics[width=1em,height=1em]{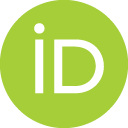} {\normalfont https://orcid.org/0000-0002-5466-0579}}
}}
\email{nfagella@ub.edu}

\address{\noindent Dept. de Matemàtiques i Informàtica\\ Universitat de Barcelona\\ Catalonia\\ Spain\\
\newline  Centre de Recerca Matemàtica\\ Bellaterra\\ Catalonia\\ Spain.
	\textsc{\newline \indent 
		\href{https://orcid.org/0000-0003-4039-5556%
		}{\includegraphics[width=1em,height=1em]{orcid2.png} {\normalfont https://orcid.org/0000-0003-4039-5556}}
}}
\email{lpardosimon@ub.edu}

\thanks{This work was partially supported by the project PID2023-147252NB-I00 financed by MICIU/AEI MCIN/AEI/10.13039/501100011033,
	FEDER, EU, by the Spanish State Research Agency, through the
	Severo Ochoa and María de Maeztu Program for Centers and Units of
	Excellence in R\&D (CEX2020-001084-M), and by the Catalan government through ICREA Academia 2020. The second author is a Serra Húnter fellow.
	}

\begin{document}
\begin{abstract}
	This article surveys the impact of Eremenko and Lyubich's paper {\em ``Examples of entire functions with pathological dynamics''}, published in 1987 in the {\em Journal of the LMS}. Through a  clever extension and use of classical approximation theorems, the authors constructed examples exhibiting behaviours previously unseen in holomorphic dynamics. Their work laid foundational techniques and posed questions that have since guided a good part of the development of transcendental dynamics.
\end{abstract}
\maketitle

\section{Introduction and historical context}
Let $f$ be a holomorphic endomorphism of $S$, where $S$ might be either the Riemann sphere, $\widehat{\C}=\C\cup \{\infty\}$,  and so $f$ is a rational map, or $S=\C,$ $f$ does not extend to $\Chat$ (i.e. infinity is an essential singularity) and $f$ is a transcendental entire map. A central theme in complex dynamics is the study of the asymptotic behaviour of the iterates, or self-compositions, $f^n$, as $n \to \infty$. Building on early foundational work by Schröder, Koenigs, Leau, Montel and others, 
the systematic study of one-dimensional holomorphic dynamics began in the early 20\textsuperscript{th} century with the pioneering contributions of Fatou and Julia. One of their key insights was the partition of the dynamical space $S$ into two complementary sets: the set where the sequence $(f^n)$ forms a normal family and the behaviour is stable under perturbation, today known as the \textit{Fatou set} and denoted $F(f)$, and its complement $J(f)=S\setminus F(f)$, the \textit{Julia set}, where the dynamics is chaotic, and small changes in initial conditions can lead to vastly different outcomes. The fundamental properties of these two sets were already established in the early works of Fatou and Julia, but several questions were posed, some of which took decades to be solved, or still remain unsolved.

Following these foundational developments, the field entered a relatively quiet period. Nevertheless, important contributions continued to appear, notably from Baker and Brolin.  A significant revival occurred in the 1980s, largely spurred by computer-generated images of the Mandelbrot set and related structures. These visualisations illustrated that even simple functions, such as quadratic polynomials, could give rise to extremely complex and aesthetically striking dynamical behaviour. Perhaps even more decisive was the new perspective brought by Sullivan, drawing connections between the iteration of rational maps and quasiconformal deformations of Riemann surfaces via his celebrated proof of the nonexistence of non-preperiodic components of the Fatou set \cite{Sullivan_82}. With these new techniques, Douady and Hubbard introduced a systematic framework for studying the dynamics of quadratic polynomials and developed powerful tools such as polynomial-like mappings \cite{Orsaynotes, polike}, leading to renormalization theory, a cornerstone  of today's modern complex dynamics. The introduction of such tools in the field and this renewed interest led to a period of intense activity and development, which continues to this day.

A \textit{singular value} of a holomorphic map $f$ is a point $w$ around which some single-valued analytic inverse branch of $f$ can be continued; equivalently, $w$ is a critical or an asymptotic value. The \textit{singular set} $S(f)$ is the closure of the set of singular values. Rational maps exemplify finiteness: they have no asymptotic values and only finitely many critical values; thus $S(f)$ is finite. Transcendental entire functions are different: infinity is an essential singularity, the (topological) degree is infinite, and $S(f)$ may be infinite; asymptotic values can occur (e.g., $0$ for $e^{z}$). Additionally, by Picard’s Theorem, in any punctured neighbourhood of infinity at most one value is omitted, hence the global behaviour is significantly more intricate and harder to control.


It is in this setting that Alexandre E. Eremenko and Mikhail Yu. Lyubich, building on Baker and Herman’s earlier work, began a systematic study of the iteration of transcendental entire functions. In a series of papers beginning in the 1980s they laid down the foundational tools and ideas for the field, e.g., see \cite{Eremenko_L_survey_1989}. Their work introduced new techniques from classical complex analysis into holomorphic dynamics and proposed a number of challenging questions that have since guided the development of the subject. Central to this narrative is the \emph{pathological examples} paper, published in English in 1987 in the \emph{Journal of the London Mathematical Society} \cite{EL_87}; its contents trace back to work undertaken in Kharkiv (Ukraine) in the fall of 1983 \cite{Erem_L_pathological_Ukraine}, with an announcement in 1984 \cite{Erem_L_Russian_1984}. Despite its modest length of 11 pages, it  had a lasting impact on the field. In this article, the authors present a refined version of Runge’s classical approximation theorem and combine it with Arakelyan’s theorem to construct five transcendental entire functions, each exhibiting surprising and at the time previously unseen dynamical behaviours.

The significance of their article lies in two aspects. First, it initiated the use of classical approximation theorems from complex analysis, particularly Runge’s and Arakelyan’s theorems, as a method for constructing transcendental entire functions with prescribed dynamical features. This approach has since become a standard tool in the field. Second, the examples provided in the paper revealed new dynamical behaviours,  considered at the time \textit{pathological} in comparison to the relatively tame and better understood world of rational dynamics, and posed natural questions, some of which still remain open. Nearly four decades later, these examples are no longer viewed as anomalies but rather as foundational and even emblematic of the richness and depth of transcendental dynamics. What once seemed exotic, unnatural or even pathological is now regarded as insightful, paradigmatic and central to our understanding of the subject.

This expository note is devoted to a careful examination of the 1987 paper, its context, techniques, and enduring influence on the modern theory of transcendental dynamics.

\section{Approximation results: the technique}
Approximation theory in complex analysis deals with how well functions that are analytic or continuous on a set can be approximated by analytic functions defined on a larger set. Classical results like Runge’s theorem (1885) show that any function analytic on a neighbourhood of a compact set $K$ with connected complement can be uniformly approximated on $K$ by polynomials. This idea was later extended by Mergelyan and Arakelyan, whose theorems allow approximation even on certain unbounded sets, as long as the sets satisfy certain regularity conditions at infinity. These results can be used to construct analytic maps with specific properties by carefully choosing the functions to be approximated. 
For more information on the topic, see the recent survey \cite{Fornaess_approx_2020}. 

In the context of dynamical systems, the underlying idea is rather simple: one wants to know whether certain model dynamics can occur for entire functions. To explore this, one prescribes the desired behaviour using a \textit{model map}, defined and being analytic in pieces (or patches) of the plane, and whose analytic extension to the whole plane might not be possible. 
The goal is to ``upgrade'' this partially defined map to an entire function, accepting some approximation error in the process. The main challenge lies in designing model maps that are flexible enough to capture the desired dynamics, while still allowing control over the approximation errors, so that the resulting entire function retains, even after iteration,  the intended properties.

Precisely to deal with the iterative aspect, one typically deals with infinitely many sets, often arising from iterated images of a given one, forcing the  errors to decease as the moduli of points get larger. This reveals a limitation of Runge’s theorem -- it only allows approximation on a single region -- as well as of Arakelyan's. To overcome this, Eremenko and Lyubich extended Runge’s result in the following way, adding some useful additional control in the process.

\begin{lemma} \label{lem:approx} {\cite[Main lemma, p. 460]{EL_87}} Let $(G_k)$ be a sequence of compact subsets of $\mathbb{C}$ with the following properties:
	\begin{itemize}
		\item[(i)] $\mathbb{C} \setminus G_k$ is connected, for every $k$;
		\item[(ii)] $G_k \cap G_m = \emptyset$, for $k \ne m$;
		\item[(iii)] $\min\{ |z| : z \in G_k \} \to \infty$ as $k\to \infty$.
	\end{itemize}
 	Let $z_k\in G_k$, $\epsilon_k>0$ and $\psi$ holomorphic on $G = \bigcup_{k=1}^\infty G_k$. Then there exists an entire function $f$ satisfying
	\begin{align*}
		|f(z) - \psi(z)| &< \epsilon_k, \quad \text{for } z \in G_k;  \\
		f(z_{k}) = \psi(z_{k}), \quad f'(z_{k}) &= \psi'(z_{k}),  \quad k\in \N.
	\end{align*}
\end{lemma}

The proof consists of successive applications of Runge's theorem: For every $m\geq 1$, they obtain a polynomial $f_m$ approximating $\psi-\sum_{i=1}^{m-1} f_i$ on $G_m$ and $0$ on $G_1\cup\ldots\cup G_{m-1} \cup K_m$, where $K_m$ is an increasing sequence of compact sets exhausting the whole plane. The errors are small enough to ensure that $\sum f_m$ converges uniformly on compact subsets of $\C$ to an entire map $f$ which also approximates $\psi$ on each $G_m$ as required. By using a geometric lemma  (suggested by Ju. I. Ljubich), the polynomials can also be chosen to satisfy
\[
\sum_{i=1}^m f_i(z_k)=\psi(z_k), \text{\ \ \ and\ \ \ } \sum_{i=1}^m f_i'(z_k)=\psi'(z_k), \quad 1\leq k\leq m.
\]
Hence $f$ satisfies the prescribed values and derivatives at the points $z_k$. 

In the coming sections we shall illustrate how this result is used to construct some of the examples. Let us note, however, that although this method, when it works, guarantees the existence of an entire function with the desired dynamics, it provides no information on the function outside the region where the model map was defined. In particular, it neither describes how the function behaves globally, nor provides any information about the number or location of  its singular values, crucial in determining the global dynamics of the map. To overcome this problem, often at the expense of flexibility, other approaches have been developed that offer more control over the resulting function. These include (dynamical) quasiconformal  surgery, see e.g. \cite{BF}, quasiconformal folding methods \cite{Bishop_2015}, constructions via Cauchy integrals \cite{Rempe_Gillen_2014,Rottenfusser_2011}, and techniques based on Hörmander’s solution to the  $\overline{\partial}$-equation \cite{Evdoridou_fast_2023}. See also \cite{Gauthier_2005} for an extension of Lemma~\ref{lem:approx} to Stein manifolds, and \cite{Bishop_2024} for a recent improvement of Runge’s theorem that allows construction of maps with a bounded set of singular values. 

\section{Examples 2 and 3: Univalent wandering and Baker domains}

\subsection{Background}
Recall that for an entire map $f$, its Fatou set $F(f)$ is the set of stability and it is defined as the domain of normality of the the family of iterates. In particular, $F(f)$ is, by definition, an open set, not necessarily connected. We call each of its connected components a \textit{Fatou component}. A Fatou component is \textit{periodic} of period $p\geq 1$ if $f^p(U)\subseteq U$, and \textit{preperiodic} if there exists $q\in \N$ such that $f^q(U)$ is contained in a periodic Fatou component. A periodic Fatou component of period one is called \textit{invariant}, and a Fatou component that is not preperiodic is called a \textit{wandering domain}. Equivalently, a wandering domain $U$ is a Fatou component with the property that $f^n(U) \cap f^m(U) \ne \emptyset$ only when $n = m$. 

The behaviour of an entire map on a periodic Fatou component is fairly well-understood. Invariant Fatou components can be classified into four types: attracting basins (when the iterates converge to an attracting fixed point) in $U$, parabolic basins (when orbits converge to a fixed point with derivative $1$ in the boundary of $U$), Siegel discs (when $U$ is conjugate via a biholomorphic map to an irrational rotation of the unit disc) and Baker domains (when $f^n\vert_{U}\to \infty$ as $n\to \infty$ and $\infty$ is an essential singularity). See, e.g. \cite[Theorem 6]{Bergweiler_1993} for details.

Observe that, by definition, Baker domains can only occur for transcendental functions. The first example was given by Fatou \cite[\S15]{Fatou}, who showed that the map $z \mapsto e^{-z}+z+1$ has a Baker domain containing the right half-plane. Similarly, wandering domains also occur exclusively for transcendental maps: Sullivan \cite{Sullivan_noWD_85} famously proved in 1985 that no rational map (and therefore no polynomial) may have a wandering domain, after Baker had already constructed the first example of a transcendental entire map with a wandering domain in 1976, \cite{Baker_wandering_76}. As a result of the absence of these types of Fatou components in the rational world, Baker and wandering domains flew under the radar for a long time, with fewer examples and a poorer understanding of their general properties.

One question of particular interest is whether the iterates on these components can be univalent or, on the contrary,  whether the degree must always be larger than one. We say that a Fatou component is \textit{univalent} if the restrictions of the iterates to this component are all univalent. A Siegel disc is, by definition, univalent, while attracting and parabolic basins cannot be univalent, as they must contain at least one singular value.  No general results of this kind had been established by the 1980s for Baker or wandering domains, and only a few examples were known. Fatou’s example with a Baker domain is not univalent, as it contains critical points. Baker’s wandering domain from 1976 is multiply connected and the degrees of the iterates increase, as shown in \cite{Baker_wandering84}. Herman provided a couple of additional examples at the time. Namely, the map $z \mapsto z + \lambda \sin (2\pi z) + 1$ is shown in \cite[\S11]{Herman_84} to have a simply connected wandering domain, for certain values of $\lambda \in \C$. Another example is $z \mapsto z + 2\pi i - 1 + e^{-z}$, described in \cite[p. 564]{Baker_wandering84}. 

These examples naturally lead to the question of whether univalent Baker and wandering domains exist. The existence of univalent wandering domains is also relevant in view of Sullivan’s proof of his no wandering domains theorem for rational maps: one first reduces to a simply connected wandering domain on which the dynamics is eventually univalent, and then propagates a Beltrami form supported there to obtain a nontrivial family of quasiconformal deformations. For rational maps this contradicts the finite dimensionality of the parameter space, so wandering domains cannot exist. For transcendental entire functions, by contrast, the space of quasiconformal deformations can be infinite dimensional; thus the existence of a wandering domain would not contradict deformation theory. Rather, it would show that Sullivan’s starting point—eventual univalence—can still occur, while the dimensional obstruction disappears.

\subsection{The examples}
Examples 2 and 3 in the paper we are studying give a positive answer to the question above: they exhibit entire functions with a univalent wandering domain and a univalent Baker domain, respectively. A sketch of the constructions is as follows. 

\begin{figure}[h]
	\centering
\includegraphics[width=\textwidth]{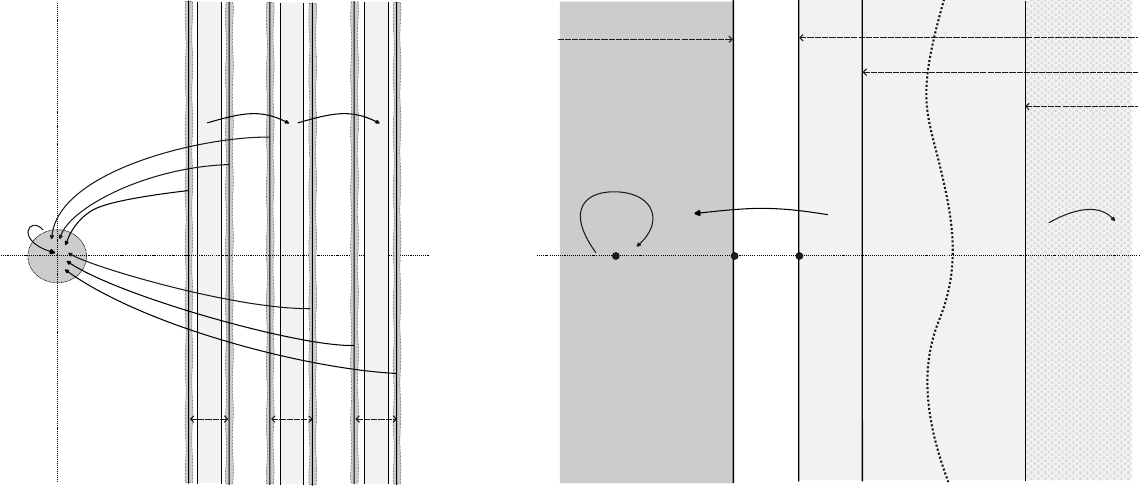}
\setlength{\unitlength}{\textwidth}
\put(-0.99,0.23){\small $\D$}
\put(-0.825,0.4){\small $\Pi^*$}
\put(-0.835,0.04){\scriptsize $\Pi_{10}$}
\put(-0.76,0.04){\scriptsize $\Pi_{11}$}
\put(-0.685,0.04){\scriptsize $\Pi_{12}$}
\put(-0.47,0.4){$e^z-6$}
\put(-0.45,0.37){$P_2$}
\put(-0.21,0.4){$2z$}
\put(-0.28,0.375){$P_1$}
\put(-0.23,0.34){$D$}
\put(-0.242,0.09){$L$}
\put(-0.08,0.31){$D_0$}
\put(-0.2,0.2){$\partial \mathcal{D}$}
\put(-0.299,0.19){\scriptsize$-4$}
\put(-0.39,0.19){\scriptsize$-5$}

\caption{\small Left: Schematic of Example~2 (univalent wandering domain): disjoint vertical strips iterate to $+\infty$; the separating vertical lines belong to an attracting basin containing the unit disc.
	Right: Schematic of Example~3 (univalent Baker domain): dynamics modelled by $z\mapsto 2z$ on a right half-plane bounded by a line $L$ that maps to a basin of attraction.%
}

\label{Figure: univalent constructions}
\end{figure}
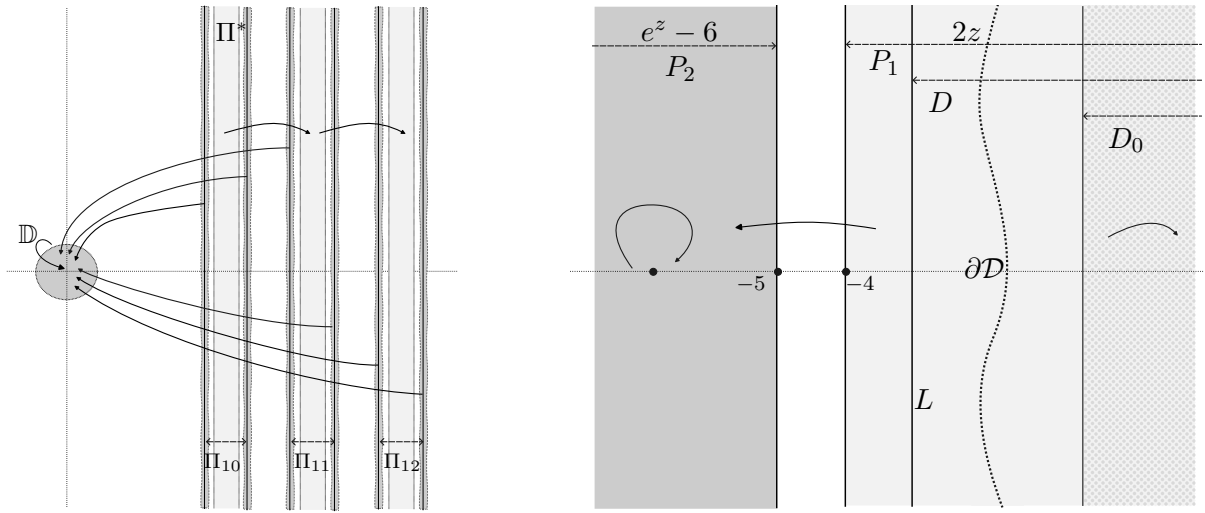

We start with Example 3 of a univalent Baker domain, the simplest construction in the paper. Our goal is to approximate the following dynamics: the univalent map $z\mapsto 2z$ is prescribed on a right half-plane, $P_1$, while attracting dynamics generated by the map $z\mapsto e^z-6$ are prescribed on a left half-plane, $P_2$ (see Figure \ref{Figure: univalent constructions}). The approximation theorem to use is Arakelyan’s theorem, since the domains are unbounded, yielding an entire function $f$ which approximates this model with an error less than $1/2$. This is enough to see that a right half-plane $D_0\subset P_1$ has ${\rm Re}(f(z))>{\rm Re}(z)+1/2$ for all $z\in D_0$, hence all orbits in $D_0$ converge to infinity and $D_0$ belongs to a Baker domain.

Injectivity of $f$ is guaranteed by the Cauchy integral formula for derivatives in a right half-plane $D\supset D_0$. The straight line $L=\partial D$ is then shown to map into $P_2$, which ensures the existence of a univalent Baker domain $\mathcal{D}$ contained in $D$ and containing $D_0$.

We point out that Herman, \cite[III, \S 9]{Herman_1985}, independently and at around the same time, gave explicit examples of functions with univalent Baker domains, namely $z \mapsto z + \frac{a}{2\pi} \sin(2\pi z) + b$ and $z \mapsto z + 2\pi i \alpha + e^z$, for suitable choices of $a, b, \alpha \in \mathbb{R}$. This fact is also mentioned in \cite{EL_87}. In the same work, Herman states additionally that if $f$ has a univalent Baker domain, then $f$ has an infinite dimensional space of quasiconformal deformations.

In the case of wandering domains, Example 2, the model map one aims to realize is defined on a sequence of disjoint closed vertical strips $\Pi''_m$, $m\geq 10$, of width increasing to 1 and centred on the line $\{z\colon \text{Re} (z) = 3m\}$. The model map is linear and expanding on each strip $\Pi''_m$ and sends it onto a vertical strip $\Pi_{m+1} \supset \Pi{_m+1}''$ of width exactly one, see Figure~\ref{Figure: univalent constructions}. One can check that the iterates of a substrip $\Pi^*\subset \Pi_{10}''$ of width $1/2$ remain inside the union of strip regions, where the model map is defined, and hence tend to $+\infty$ along the real direction.

A subtlety here is that, even if one succeeds in extending such a model map to an entire function, there is no guarantee that the resulting strips will lie in distinct Fatou components: they could all be part of a single Baker domain. To prevent this, the model map also sends the closed unit disc, $\overline{\D}$, compactly inside itself (ensuring the existence of an attracting fixed point inside) and the boundaries of the strips $\Pi_m$ inside $\D$. This forces the wandering vertical strips to lie in different Fatou components (of the final entire map), as desired. See Figure \ref{Figure: univalent constructions}. This key idea is further refined and exploited in subsequent publications, as we discuss in Section 4. 

The approximation theorem used in this case is also Arakelyan’s theorem (with a slight modification to allow errors $\epsilon_m \to 0$ on the strips $\Pi_m''$, as in Lemma \ref{lem:approx}), since the regions where the model is defined are not compact. 

The transfer of univalence from the model map to the resulting entire function is made possible by carefully controlling the approximation errors—chosen sufficiently small on each $\Pi_m''$— and using the Cauchy integral formula for derivatives to ensure that the derivative does not vanish and the distortion remains within the required bounds. The entire function obtained therefore has a wandering domain $\mathcal{D}\supset \Pi^*$.

\subsection{Subsequent developments}
 Eremenko and Lyubich's approach was later used to construct further examples of Baker domains. In \cite{Rippon_1999_2}, univalent Baker domains of arbitrary period $p \geq 1$ are constructed, while \cite{Baranski_2001} provides examples of \textit{spiralling} Baker domains.  Incidentally,  the only unbounded Fatou components whose boundary may be a Jordan curve in $\widehat{\C}$ are precisely univalent Baker domains \cite{Baker_Weinreich_1991}. 
 
The understanding of Baker domains—both univalent and not—has advanced considerably since the example of Eremenko and Lyubich.  Building on Cowen’s results on the iteration of holomorphic functions on the unit disc \cite{Cowen}, a complete classification of simply connected Baker domains has been obtained, together with a description of the space of quasiconformal deformations supported on this type of stable component, see \cite{Koning_99, Fagella_2006}. Further related examples were later constructed in \cite{Kisaka_1998, Fagella_2006, Bergweiler_2012}, and a comprehensive overview of later developments concerning Baker domains can be found in  Rippon's survey \cite{Rippon_2008_survey}.

Likewise,  wandering domains are  generally not univalent. In fact, multiply connected wandering domains are never univalent; see \cite{Bergweiler_PG_multiply_13} for a description of the dynamics in multiply connected wandering domains. Simply connected wandering domains, on the other hand, exhibit a much wider range of behaviours, as established in \cite{Benini_2022}. We shall discuss these results in the following section. Further examples of wandering domains for meromorphic maps constructed following this technique include \cite{Singh_2000, Baker_90, Rippon_2008}. For analogous constructions to those in \cite{EL_87} yielding self-maps of $\C^*$, see  \cite{Kotus_90}.

\section{Example 1: an oscillating wandering domain}

The examples of wandering domains discussed so far are \emph{escaping}, in the sense that all their iterates tend to infinity. Formally, this means that the family of iterates, $(f^n)_{n\geq 0}$, when restricted to the wandering domain, converges to the constant function taking the value infinity. In fact, it follows from standard arguments that all limit functions  (i.e. limits of subsequences $(f^{n_k})_k$)  in a wandering domain are constant \cite{Fatou_1920}. However, there is {\em a priori} no reason to expect that all such subsequences converge to infinity -- the pairwise disjoint Fatou components in the orbit of the wandering domain could have finite limit functions as well, either in addition to or instead of infinity.

\begin{figure}[htb!]
\centerline{
\includegraphics[width=0.7\textwidth]{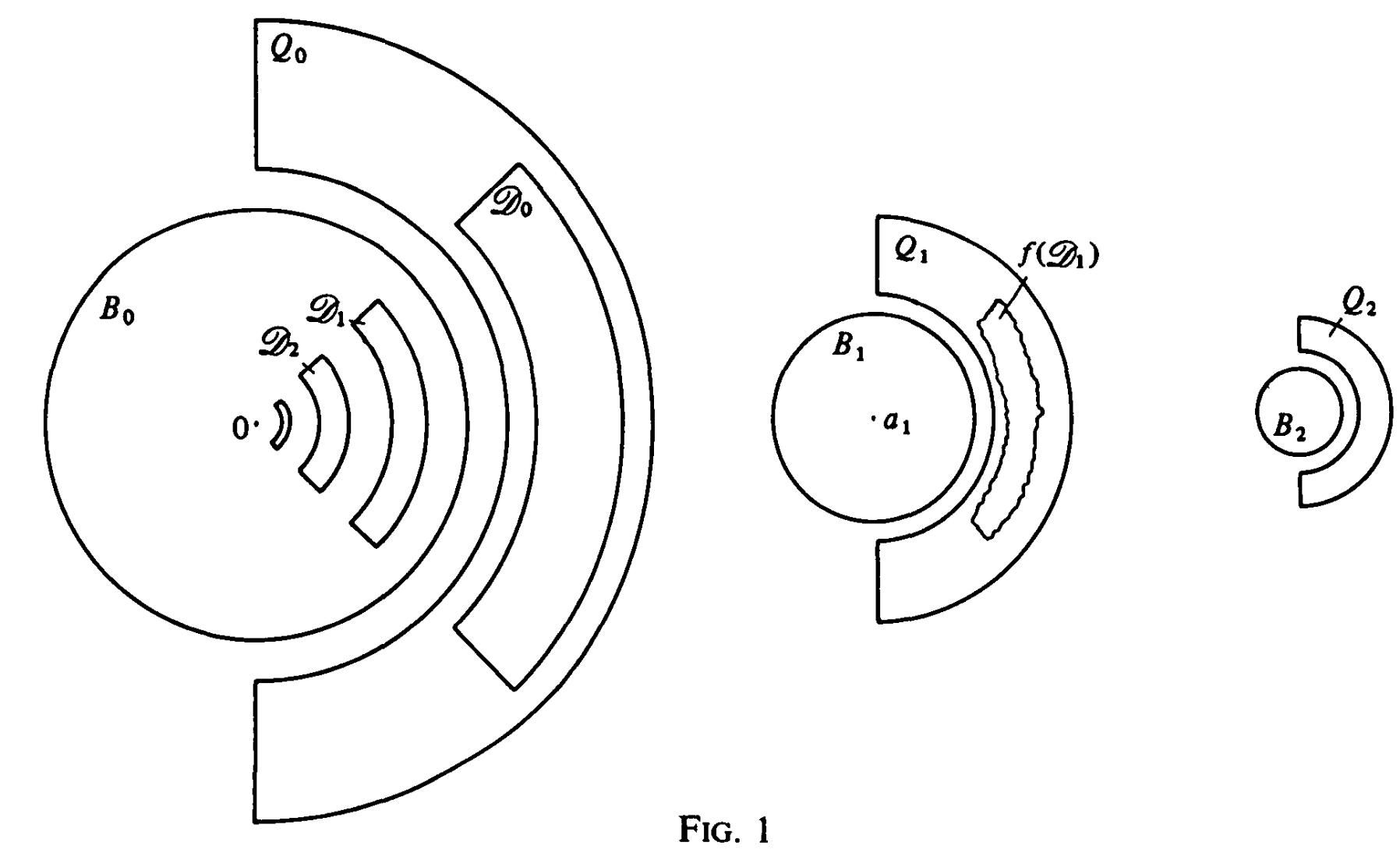}
}
\caption{\label{Fig:oscillating_orig} \small  The original figure in \cite{EL_87} illustrating the construction of an oscillating wandering domain containing $\mathscr{D}_0$. For $m\geq 0$, the model map translates every ball $B_m$ to the right so that $a_m$ is sent to $a_{m+1}$.  After $m$ iterates, $\mathscr{D}_m$ maps inside $Q_m$, which in turn is mapped inside $\mathscr{D}_{m+1}$, creating an oscillating pattern. Sets and errors are carefully chosen so that the resulting approximation map satisfies the same inclusions, and hence $a_0=0$ is a limit point for every orbit in $\mathscr{D}_0$. Hence each $a_m$ and $\infty$ occur as limit functions. }
\end{figure}

Example 1 in \cite{EL_87} provides the first (and still one of the few known) instances of a wandering domain with infinitely many limit functions,  infinity being one of them. In modern terminology, such a domain is known as an \emph{oscillating} wandering domain, meaning that the iterates accumulate both at infinity and at least one (and thus infinitely many) other finite value. The rough idea of their construction is as follows: using their ``infinite version'' of Runge's theorem, Lemma \ref{lem:approx}, the authors prescribe translations between a sequence of balls $(B_m)$ tending to infinity, with radii $r_m\to 0$; see Figure \ref{Fig:oscillating_orig}. Around each of them, a half-annulus $Q_m$  is mapped to an appropriate domain $\mathscr{D}_{m+1}$ in $B_0$, successively closer to $0$, ensuring that iterates of a neighbourhood of $\mathscr{D}_m$ remain within the collection of balls for sufficiently long, and that $f^{m+1}(\mathscr{D}_{m+1})$ is contained in the half-annulus $Q_{m+1}$, creating an  \textit{oscillating} behaviour. All parameters in the construction, including radii and errors, are chosen with sufficient care to ensure that the resulting entire map has a Fatou component, necessarily a wandering domain, replicating the prescribed dynamics.

Following their construction, the authors posed the question of whether there exists an entire function with a wandering domain for which all limit functions are finite. This question remains open to this day and is considered one of the most prominent open problems in transcendental dynamics. Various formulations and variants of this question, attributed to the authors and others, along with progress towards their resolution, are compiled in Hayman and Lingham's open problems list \cite[Problems 2.67, 2.77, and 2.87]{Hayman_book_50}.

\subsection{Subsequent developments}
As mentioned earlier, the construction methods in \cite{EL_87} give no information about the global features of the resulting map $f$. In particular, one knows neither the number nor the location of its singular values. The singular set $S(f)$ plays a crucial role in holomorphic dynamics; indeed, every stable component is associated with at least one singular value, so one may say that this set of ``special'' points governs the dynamics of $f$; see \cite[\S 4.3]{Bergweiler_1993}, \cite{Bergweiler_1995} or \cite{Baranski_2020}.

For polynomial and rational maps, this set is always finite. It is therefore natural to look for parallels between the dynamics of such maps and those in the \textit{Speiser class} $\mathcal{S}$, consisting of all transcendental entire functions with a finite number of singular values. This class, which includes exponential and trigonometric maps, retains some of the ``tameness'' characteristics of rational dynamics. For example, Eremenko and Lyubich \cite{Erem_L_Russian_1984} (see also \cite[Theorem~3]{Erem_L_classes_1992}), and independently Goldberg and Keen \cite[Theorem 4.2]{Goldberg_Keen_86}, proved in the 1980s that functions in $\mathcal{S}$ do not have wandering domains. A natural extension of the class $\mathcal{S}$ is the \textit{Eremenko–Lyubich class} $\mathcal{B}$, consisting of all entire functions for which $S(f)$ is bounded. Shortly after, Eremenko and Lyubich showed that functions in $\mathcal{B}$ cannot have escaping wandering domains \cite[Theorem 1]{Erem_L_classes_1992}, although it remained an open question for nearly twenty years whether they could have wandering domains at all.

From the 1990s onwards, building on the foundational work of Eremenko and Lyubich, research in transcendental dynamics largely focused on specific families of functions, such as exponential and cosine maps, or on general properties of functions in the class $\mathcal{B}$. For overviews of the dynamics in this class, see \cite{Sixsmith_survey_2018} and the recent survey on escaping dynamics \cite{Bergweiler_Rempe_escaping_survey_25}. During this period, the study of wandering domains continued but received comparatively less attention until 2015, when C. Bishop developed his technique of quasiconformal folding and used it, among other things, to construct the first example of a function in $\mathcal{B}$ with an oscillating wandering domain \cite{Bishop_2015}. This result revived interest in the topic, and several variants of Bishop’s construction soon followed, yielding functions in $\mathcal{B}$ with additional features such as unbounded oscillating wandering domains \cite{Lazebnik_2017} and univalent ones, \cite{Fagella_2019}. Subsequently, using quasiregular interpolation, Martí-Pete and Shishikura \cite{Marti_Pete_2020} constructed functions in $\mathcal{B}$ with oscillating wandering domains with the additional property of being of finite order of growth. In parallel, Bishop also constructed, via an infinite product, a function with multiply connected wandering domains whose Julia set has Hausdorff dimension 1, \cite{Bishop_dim1_18}, thereby solving a long-standing problem in transcendental dynamics.

These developments renewed momentum in the study of wandering domains and brought fresh attention to the techniques introduced by Eremenko and Lyubich in the paper at hand, which have recently found striking and unexpected applications that we now describe.

The classification of periodic Fatou components outlined in Section~3 dates back to Fatou’s work a century ago. However, certain Fatou components admit a finer classification. This is the case of Baker domains  \cite{Fagella_2006} and much more recently, of simply connected wandering domains \cite{Benini_2022} and multiply connected ones \cite{Gustavo}. Benini, Evdoridou, Fagella, Rippon, and Stallard \cite{Benini_2022} provide a framework for understanding the different types of behaviour that can occur in these domains. Their approach takes into account both the hyperbolic distances between iterates and how orbits interact with the boundaries of the wandering domains, leading to a classification into nine distinct dynamical types. The construction of examples realizing each of the nine different types  is heavily inspired by the techniques developed in \cite{EL_87}. In particular, they use a slightly refined version of  Lemma \ref{lem:approx} that allows for the prescription of additional orbits. Using this extended version of the lemma, the authors construct escaping wandering domains that are, roughly speaking, translated copies of a disc, with extra care taken to control the location of the boundaries, proving that these could even be Jordan curves.  The exact prescription of a certain number of points allowed by the key Lemma \ref{lem:approx} plays a crucial role here, since it determines the type of wandering domain that is obtained.

Further constructions using techniques from approximation theory include examples of bounded oscillating wandering domains \cite{Evdoridou_2023}, and of unbounded wandering domains exhibiting different internal dynamics \cite{Evdoridou_fast_2023}.

The control of boundaries in the construction of wandering domains was taken significantly further by Boc Thaler in \cite{BocThaler_2021}, who showed that any bounded, connected, regular open set whose closure has connected complement can be realized as an escaping (or as an oscillating) univalent wandering domain of an entire function. Recall that a domain is regular if it coincides with the interior of its closure. In particular, this result implies that the boundaries of bounded simply-connected wandering domains can be arbitrarily well-behaved, for instance, even analytic, in contrast to invariant Fatou components, whose boundaries are known to be generally highly non-smooth; see \cite{Azarina_89} and \cite[Theorem C]{Baranski_2025}.

Perhaps rather surprisingly, the exact prescription of the boundary is again achieved using the approximation techniques introduced by Eremenko and Lyubich. However, Boc Thaler introduces a crucial new idea: rather than approximating the entire forward orbit of a domain in a single step using Lemma \ref{lem:approx}, the construction proceeds inductively by applying  some version of the lemma infinitely many times. More precisely, a sequence $(f_n)$ of entire functions is constructed, converging uniformly on compact sets to a limit map $f$. At each step $n$, the map $f_n$ is defined so as to approximate $f_{n-1}$ on a large disc that contains the images of the wandering domain under the first $n-1$ iterates. Additionally, $f_n$ is carefully defined on a new region where the $n$-th iterate will lie, so as to shape the domain in the desired way by mapping ``layers'' of points to an attracting basin. This inductive control allows the construction to ``sculpt'' the boundaries of the resulting domain precisely.  

Building on Boc Thaler’s techniques, Martí-Pete, Rempe and Waterman  \cite{MartiPete_JAMS_2025}  pushed the boundary control even further, showing that any bounded simply connected domain $D$ such that $\C\setminus \overline{D}$ has an unbounded connected component $W$ with $\partial W = \partial D$ can be realized as an escaping wandering domain of a transcendental entire function. This significantly broadens the class of admissible domains, allowing for examples with highly intricate geometry, like for example  a \textit{Lakes of Wada continuum}, where the complement $\C \setminus \overline{D}$ has two or more (possibly infinitely many) connected components, all of which share the same boundary $\partial D$. This marks the first appearance of Lakes of Wada continua in holomorphic dynamics. For related results in the meromorphic setting, see \cite{MartiPete_merom_2025}.

The constructions in \cite{BocThaler_2021, MartiPete_JAMS_2025} are further modified in \cite{Pardo_2024} to produce an example of an oscillating wandering domain with the
property that, in a precise sense, nearly all of its forward iterates are contained within a bounded domain, a problem closely related to the fundamental question of existence of  wandering domains with all its limit functions in a bounded set.

But perhaps the most unexpected and celebrated result building on the paper we are surveying is the construction, in the same article \cite{MartiPete_JAMS_2025}  by Martí-Pete, Rempe and Waterman, of a counterexample to \textit{Eremenko's conjecture}. This conjecture, posed by Eremenko in 1989 \cite{Eremenko_1989}, asserts that every connected component of the \textit{escaping set} of a transcendental entire function, defined as 
\begin{equation}\label{eq_If}
	I(f)=\left\{ z\in \C\colon f^n(z)\to \infty \text{ as }n\to \infty\right\},
\end{equation}
is unbounded. This conjecture has been a driving force in the field for many years, motivating a large body of results and establishing connections to other important questions; see the introduction of \cite{MartiPete_JAMS_2025} for details.

\section{Examples 4 and 5: Julia sets of positive area supporting invariant line fields}

For the last two examples constructed in \cite{EL_87}, attention shifts from the stability of the Fatou set to the dynamics on its complement, the Julia set. Due to the expansive nature of a meromorphic function $f$ on its Julia sets, $J(f)$ is either the entire complex plane (or, if $f$ is rational, the Riemann sphere $\widehat{\mathbb{C}}$), or it has empty interior. Examples of rational functions with $J(f) = \widehat{\mathbb{C}}$ arise from the multiplication theorems of elliptic functions, particularly the class of Lattès maps, see, e.g. \cite[Theorem 7.3]{Milnor_book}. For entire functions, the first known example was given by Baker, \cite{Baker_limit_70}, who showed that $J(\lambda z e^z) = \mathbb{C}$ for a suitable value of $\lambda$. Later, Misiurewicz \cite{Misiurewicz_81} confirmed a conjecture of Fatou by proving that $J(e^z) = \mathbb{C}$.

\begin{figure}[h]
	\centering
	\def\svgwidth{0.95\linewidth}
\begingroup%
  \makeatletter%
  \providecommand\color[2][]{%
    \errmessage{(Inkscape) Color is used for the text in Inkscape, but the package 'color.sty' is not loaded}%
    \renewcommand\color[2][]{}%
  }%
  \providecommand\transparent[1]{%
    \errmessage{(Inkscape) Transparency is used (non-zero) for the text in Inkscape, but the package 'transparent.sty' is not loaded}%
    \renewcommand\transparent[1]{}%
  }%
  \providecommand\rotatebox[2]{#2}%
  \newcommand*\fsize{\dimexpr\f@size pt\relax}%
  \newcommand*\lineheight[1]{\fontsize{\fsize}{#1\fsize}\selectfont}%
  \ifx\svgwidth\undefined%
    \setlength{\unitlength}{528bp}%
    \ifx\svgscale\undefined%
      \relax%
    \else%
      \setlength{\unitlength}{\unitlength * \real{\svgscale}}%
    \fi%
  \else%
    \setlength{\unitlength}{\svgwidth}%
  \fi%
  \global\let\svgwidth\undefined%
  \global\let\svgscale\undefined%
  \makeatother%
  \begin{picture}(1,0.31107955)%
    \lineheight{1}%
    \setlength\tabcolsep{0pt}%
    \put(0,0){\includegraphics[width=\unitlength,page=1]{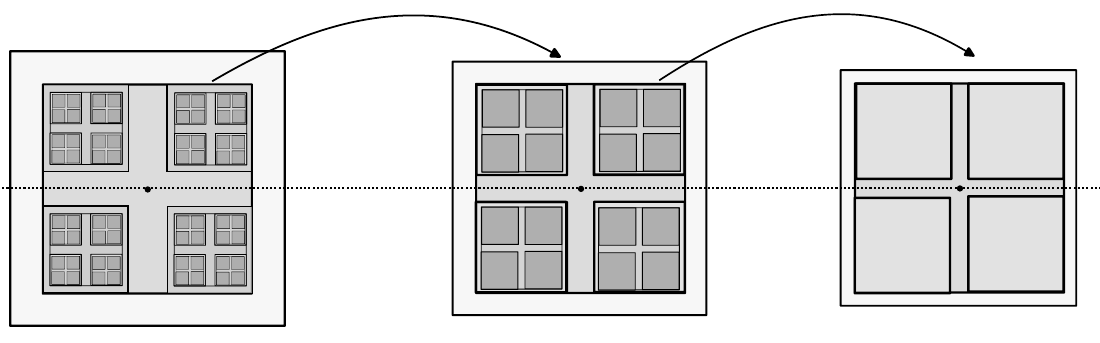}}%
    \put(0.12280686,0.14679377){\color[rgb]{0,0,0}\makebox(0,0)[lt]{\lineheight{1.25}\smash{\begin{tabular}[t]{l}\fontsize{9pt}{1em}$4k$\end{tabular}}}}%
    \put(0.51632072,0.26969488){\color[rgb]{0,0,0}\makebox(0,0)[lt]{\lineheight{1.25}\smash{\begin{tabular}[t]{l}\fontsize{10pt}{1em}$Q^1_{k+1}$\end{tabular}}}}%
    \put(0.79394695,0.18215466){\color[rgb]{0,0,0}\makebox(0,0)[lt]{\lineheight{1.25}\smash{\begin{tabular}[t]{l}\fontsize{10pt}{1em}$Q_{1,k+2}$\end{tabular}}}}%
    \put(0.89421013,0.18127293){\color[rgb]{0,0,0}\makebox(0,0)[lt]{\lineheight{1.25}\smash{\begin{tabular}[t]{l}\fontsize{10pt}{1em}$Q_{2,k+2}$\end{tabular}}}}%
    \put(0.78881997,0.08024645){\color[rgb]{0,0,0}\makebox(0,0)[lt]{\lineheight{1.25}\smash{\begin{tabular}[t]{l}\fontsize{10pt}{1em}$Q_{3,k+2}$\end{tabular}}}}%
    \put(0.89523376,0.08042603){\color[rgb]{0,0,0}\makebox(0,0)[lt]{\lineheight{1.25}\smash{\begin{tabular}[t]{l}\fontsize{10pt}{1em}$Q_{4,k+2}$\end{tabular}}}}%
  \end{picture}%
\endgroup%

	\caption{\small Schematic of the construction of Example 4, a Julia set of positive measure which is not the entire complex plane. }
	\label{Figure: squares}
\end{figure}
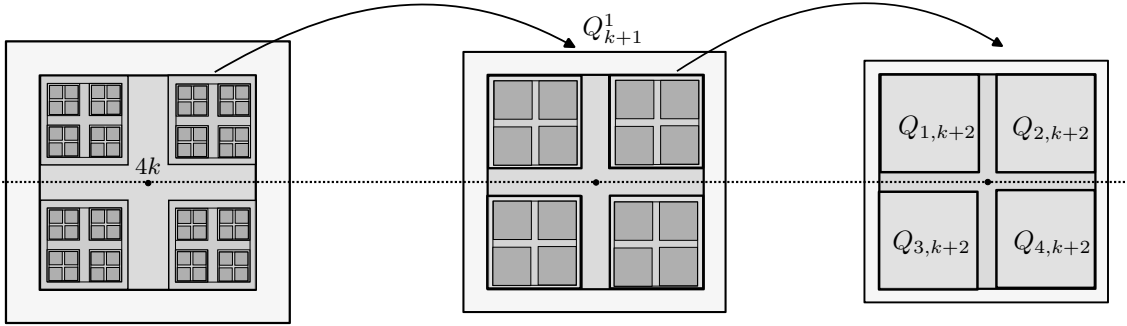

Example 4 in \cite{EL_87} provides an instance of an entire function $f$ for which $J(f) \neq \mathbb{C}$, yet $J(f)$ has positive Lebesgue measure. The idea behind the construction is as follows. Consider a sequence of pairwise disjoint unit squares $Q_k$ arranged along the real axis, each containing the point $4k$, for $k \geq 0$. Around each $Q_k$, take a slightly larger square, $Q_k^1$, of side length $1 + \epsilon_k$, and inside each $Q_k$, select four smaller squares $Q_{j,k}$ (for $j = 1, \dots, 4$), arranged so that the remaining set $E_k := Q_k \setminus \bigcup_j Q_{j,k}$ has area tending to zero as $k \to \infty$, see Figure \ref{Figure: squares}. Using Lemma \ref{lem:approx}, one constructs an entire function $f$ that maps each $Q_{j,k}$ univalently onto $Q_{k+1}$ with sufficiently large derivative to ensure strong expansion. An attracting fixed point is also introduced in the construction to guarantee that the Julia set is not the entire complex plane. 

The set of points in $Q_0$ that remain in the sequence of squares under iteration, i.e., those obtained by taking successive preimages of the $Q_{j,k}$ as $k \to \infty$, forms a Cantor set $K$. By construction, $K \subset J(f)$, since the prescribed large derivatives enforce expansion on this set. Moreover, $K$ has positive Lebesgue measure due to two key features: the strong contraction of preimages (caused by the large derivatives) and the fact that the area of the complement of the iterates of $K$ in the union of squares is controlled by the measure of the sets $E_k$, $k\geq 0$, which tends to zero as $k\to \infty$.

Observe that $K$ is not only contained in $J(f)$, but, by construction, it also lies within the escaping set $I(f)$. In parallel with the work in \cite{EL_87}, McMullen, \cite{McMullen_area_87}, showed that all maps in the sine family, given by $z \mapsto \sin(\alpha z+\beta)$, $\alpha \neq 0$, also exhibit this property: their escaping set is contained in the Julia set  and has positive Lebesgue measure.

The authors show, additionally, that the map $f$ in Example 4 supports a measurable invariant line field on $J(f)$. Following \cite{EL_87}, for a field of straight lines defined on almost all $z\in J(f)$, let $\theta(z)\in\mathbb{R}/\pi\mathbb{Z}$ be the angle of the line at $z$ with the positive real axis, and encode it by
	$$
	\mu(z)=
	\begin{cases}
		e^{2i\theta(z)}, & z\in J(f),\\[2mm]
		0, & z\notin J(f).
	\end{cases}
	$$
	The field is measurable if $\mu$ is measurable, and it is non-trivial if $\operatorname{mes}(J(f))>0$. It is $f$-invariant when
	\begin{equation}\label{eq:invariant}
	\mu\big(f(z)\big)=\left(\frac{f'(z)}{|f'(z)|}\right)^{2}\mu(z) \quad \text{ for almost all }z\in J(f).
	\end{equation}

Equivalently, an invariant line field is a Beltrami coefficient $\mu(z)$ supported on $J(f)$  such that $|\mu(z)|=1$  and $f^*\mu=\mu$. The relation between the two definitions is given by $\mu(z)=\exp(2 i \theta(z))$.

The existence of invariant line fields supported on the Julia set of an entire or meromorphic function (which, by definition, implies that $J(f)$ has positive measure) is central to the conjecture known as \textit{density of hyperbolicity}, which asserts that hyperbolic functions are dense in natural families of meromorphic functions under suitable assumptions, one of the main problems in one-dimensional real and holomorphic dynamics. For example, in the quadratic family $z\mapsto z^2+c$, $c\in \C$, the existence of non-hyperbolic components of the interior of the Mandelbrot set is equivalent to the existence of polynomials with positive-measure Julia sets that carry an invariant line field (see e.g. \cite{mss,McMullen_renormalization}). Although it was long believed that Julia sets of quadratic polynomials might always have zero area, which would rule out invariant line fields, this was disproved by Buff and Cheritat \cite{Buff_Cheritat_area_2012}. In transcendental dynamics the relationship is less direct: For a family whose set of singular values is unbounded, not even the set of structurally stable maps (a weaker condition than hyperbolicity) must be dense -- see for example \cite{ABF1}.

In the paper at hand, Eremenko and Lyubich show that the function built in Example 4 carries an invariant line field in $K$, a subset of its Julia set. By adjusting the errors of the approximation they show that 
\[
|\arg f'(z)| \leq A \, \delta_k, \quad \quad \text{for all $z\in Q_{k,j}$}, 
\]
where $A$ is a constant and $\delta_k$ is a summable sequence. This implies that the series $\theta(z):=-\sum_{k\geq 0} \arg f'(f^k(z))$ converges absolutely and uniformly on the union of forward images of the Cantor set and defines a continuous function, which satisfies the invariance formula 
\[
\theta(f(z))=\theta(z)+\arg f'(z),
\]
or equivalently, satisfying \eqref{eq:invariant}. Hence $\theta$ is a non-trivial invariant line field on a subset of the Julia set of positive measure.

Furthermore, Example 5 (a modification of Example 4) constructs a function with an infinite-dimensional family of invariant line fields, something impossible for functions in class $\mathcal{S}$.

\subsection{Subsequent developments}
Further families and classes of transcendental entire functions have been shown to have Julia sets of positive area, by proving that the escaping set is contained in the Julia set and itself has positive measure. This strategy is particularly effective when $f \in \mathcal{B}$, since for such functions it is always the case that $I(f) \subset J(f)$. In some cases, it is even shown that the complement of $I(f)$ (or $J(f)$) in unbounded regions has finite area. Along these lines, Schubert~\cite{Schubert_08} showed that the Fatou set of the sine map has finite area in vertical strips of width $2\pi$, thereby complementing McMullen’s result. Hemke~\cite{Hemke_2005}, and later Wolff~\cite{Wolff_exponentials_21}, 
studied families of functions defined as linear combinations of exponentials and polynomials, providing conditions under which $\mathrm{meas}(\mathbb{C} \setminus I(f)) < \infty$.

Sixsmith~\cite{Sixsmith_area_15} proved that for certain families generalizing exponential and sine functions, the escaping set, and even its subset, the \emph{fast} escaping set, has positive measure and exhibits a particular topological structure known as a \textit{spider's web}. Additional sufficient conditions for $I(f)$ to have positive measure, stated in terms of growth estimates for functions in class $\mathcal{B}$, were given by Aspenberg and Bergweiler~\cite{Aspenberg_12}, and were later shown to be sharp by Cui~\cite{Cui_2018}. Further results for functions outside class $\mathcal{B}$ with Julia sets of positive measure, based on growth and regularity conditions, are given in~\cite{Bergweiler_16, Bergweiler_18}.

In  McMullen and Eremenko-Lyubich examples, the Julia set has infinite positive measure. In~\cite{Wolff_2025}, Wolff refines the constructions in~\cite[Examples 4 and 5]{EL_87}, using a more intricate grid of squares with carefully prescribed dynamics to construct a transcendental entire function whose Julia set has finite positive measure. 

Functions in class $\mathcal{S}$ may have Julia sets, and even escaping sets, of positive area, as illustrated by the example $z\mapsto \sin(z)$ discussed above. Nevertheless, it is conjectured that functions in $\mathcal{S}$ admit no invariant line fields supported on their Julia set, in analogy to the so-called \textit{Fatou's conjecture } in rational dynamics, see \cite[\S8, p. 1016]{Erem_L_classes_1992}. The absence of line fields supported on Julia sets has been proved for certain classes of non-recurrent meromorphic functions in class $\mathcal{B}$ \cite{Rempe_2011}, as well as for meromorphic maps satisfying specific conditions on the postcritical set \cite{Graczyk}. Moreover, it is known that functions in class $\mathcal{B}$ cannot support invariant line fields on their escaping set, unlike in \cite[Example 4]{EL_87}; see \cite[Theorem 1.2]{Rempe_2009} and also \cite[\S 8.6]{Bergweiler_Rempe_escaping_survey_25}. Several results concerning functions whose Julia set has zero area (and therefore cannot support invariant line fields) can be found, for example, in \cite{Erem_L_classes_1992, Stallard_zero_90, Wolff_22, Jankowsi_97, Zheng_02}.

\medskip
\noindent \textbf{Acknowledgements.} We are grateful to Walter Bergweiler, Phil Rippon and the referees for directing us to relevant references and for valuable feedback that improved the exposition.

\bibliographystyle{alpha}
\bibliography{references}

\end{document}